\documentclass[final,twocolumn]{IEEEtran} 

\usepackage{epsfig}
\usepackage{amsmath}
\usepackage{amssymb} 

\title{\LARGE {\bf Steering laws for motion camouflage 
}}

\author{ \parbox{3 in}{\centering E. W. Justh \\
         Institute for Systems Research \\
         University of Maryland \\
         College Park, MD 20742, USA \\
         {\tt\small justh@umd.edu}}
         \hspace*{ 0.5 in}
         \parbox{3 in}{ \centering P. S. Krishnaprasad \\
         Institute for Systems Research and \\
         Dept. of Electrical and Computer Engineering \\
         University of Maryland \\
         College Park, MD 20742, USA \\
         {\tt\small krishna@umd.edu}}
}

\begin{document}

\maketitle

\begin{abstract}

Motion camouflage is a stealth strategy observed in nature.
We formulate the problem as a feedback system for
particles moving at constant speed, and define what 
it means for the system to be in a state of motion camouflage.
(Here we focus on the planar setting, although the results 
can be generalized to three-dimensional motion.) 
We propose a biologically plausible feedback law, and use
a high-gain limit to prove accessibility of a motion camouflage
state in finite time.  We discuss connections to work in missile
guidance. We also present simulation results to explore the 
performance of the motion camouflage feedback law for a variety
of settings.

\end{abstract}

\section{Introduction}

Motion camouflage is a stealth strategy employed by various visual
insects and animals to achieve prey capture, mating or territorial  
combat. In one type of motion camouflage, the predator 
camouflages itself against a fixed background object
so that the prey observes no relative motion between 
the predator and the fixed object.
In the other type of motion camouflage, the predator 
approaches the prey such that from the point
of view of the prey, the predator always appears to
be at the same bearing.  (In this case, we say that the
object against which the predator is camouflaged is
the point at infinity.)  Assuming that the prey can 
readily observe optical flow, but only poorly sense
looming, this type of motion by the predator is then
difficult to detect by the prey.  For example, insects
with compound eyes are quite sensitive to optical flow
(which arises from the transverse component of the relative
velocity between the predator and the prey), but are
far less sensitive to slight changes in the size of
images (which arise from the component of the
relative velocity between the predator and prey along
the line between them). More broadly such interactions may also apply
in settings of mating activity or territorial maneuvers as well. 
In the work, \cite{srinivasan95} of 
Srinivasan and Davey, it was suggested that the data on visually
mediated interactions between two hoverflies, {\em Syritta pipiens}
obtained earlier by Collett and Land \cite{collett75}, supports
a motion camouflage hypothesis. Later, Mizutani, Chahl and Srinivasan 
\cite{srinivasan03}, observing territorial aerial maneuvers of 
dragonflies {\em Hemianax papuensis}, concluded that the flight
pattern is motivated by motion camouflage (see Figure 1 in their paper).  
See also \cite{srinivasan04} for a review of related themes in
insect vision and flight control.

Motion camouflage can be used by a predator to
stealthily pursue prey, but a motion camouflage strategy
can also be used by the prey to evade a predator.
The only difference between the strategy of the predator
and the strategy of the evader is that the predator seeks
to approach the prey while maintaining motion camouflage,
whereas the evader seeks to move away from the predator
while maintaining motion camouflage.  
Besides explaining certain biological pursuit strategies,
motion camouflage may also be quite useful in certain military
scenarios (although the ``predator'' and ``prey'' labels
may not be descriptive). In some settings, as is the case 
in \cite{srinivasan95}, \cite{collett75}, \cite{srinivasan03}
it is more appropriate to substitute the labels ``shadower'' 
and ``shadowee'' for the predator-prey terminology.

In this work, we take a structured approach to deriving
feedback laws for motion camouflage, which incorporate
biologically plausible (vision) sensor measurements.
We model the predator and prey as point particles
moving at constant (but different) speeds, and subject to
steering (curvature) control.
For an appropriate choice of feedback control law for
one of the particles (as the other follows a prescribed
trajectory), a state of motion camouflage is then approached
as the system evolves.
(In the situation where the predator follows a motion-camouflage
law, and the speed of the predator exceeds the speed of the
prey, the predator is able to pass ``close'' to the prey in
finite time.  In practice, once the predator is sufficiently
close to the prey, it would change its strategy from a 
pursuit strategy to an intercept strategy.)

What distinguishes this work from earlier study of 
motion-camouflage trajectories in 
\cite{glendinning} is that we present {\it biologically
plausible feedback laws} leading to motion camouflage.
Furthermore, unlike the neural-network approach used
in \cite{anderson} to achieve motion camouflage using
biologically-plausible sensor data, our approach gives
an explicit form for the feedback law which has a 
straightforward physical interpretation.

The study of motion camouflage problems also naturally 
extends earlier work on interacting systems 
of particles, using the language of curves and 
moving frames \cite{scl02}-\cite{cdc05}.     

\section{Planar pursuit-evasion model}

\begin{figure}
\hspace{.5cm}
\epsfxsize=7cm
\epsfbox{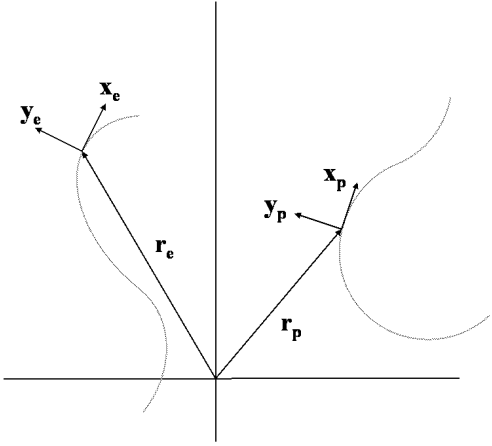}
\caption{\label{framefig2d} Planar trajectories for the pursuer
and evader, and their respective natural Frenet frames.}
\end{figure}

For concreteness, we consider the problem of motion
camouflage in which the predator (which we refer to 
as the ``pursuer'') attempts to intercept the prey
(which we refer to as the ``evader'') while appearing
to the prey as though it is always at the same bearing  
(i.e., motion camouflaged against the point at infinity).
In the model we consider, the pursuer moves at unit
speed in the plane, while the evader moves at a constant
speed $ \nu < 1 $. 
The dynamics of the pursuer are given by
\begin{eqnarray}
\label{pursuer2d}
\dot{\bf r}_p \hspace{-.2cm} & = & \hspace{-.2cm}
 {\bf x}_p, \nonumber \\
\dot{\bf x}_p \hspace{-.2cm} & = & \hspace{-.2cm}
 {\bf y}_p u_p,  \nonumber \\
\dot{\bf y}_p \hspace{-.2cm} & = & \hspace{-.2cm}
 -{\bf x}_p u_p, 
\end{eqnarray}
where $ {\bf r}_p $ is the position of the pursuer, 
$ {\bf x}_p $ is the unit tangent vector to the trajectory of
the pursuer, $ {\bf y}_p $ is the corresponding unit normal vector
(which completes a right-handed orthonormal basis with $ {\bf x}_p $),
and the plane curvature $ u_p $ is the steering control for the
pursuer.  Similarly, the dynamics of the evader are
\begin{eqnarray}
\label{evader2d}
\dot{\bf r}_e \hspace{-.2cm} & = & \hspace{-.2cm}
 \nu {\bf x}_e, \nonumber \\
\dot{\bf x}_e \hspace{-.2cm} & = & \hspace{-.2cm}
 \nu {\bf y}_e u_e, \nonumber \\
\dot{\bf y}_e \hspace{-.2cm} & = & \hspace{-.2cm}
 -\nu{\bf x}_e u_e, 
\end{eqnarray}
where $ {\bf r}_e $ is the position of the evader,
$ {\bf x}_e $ is the unit tangent vector to the trajectory of 
the evader, $ {\bf y}_e $ is the corresponding unit normal vector,
and $ u_e $ is the steering control for the evader.
Figure \ref{framefig2d} illustrates equations (\ref{pursuer2d})
and (\ref{evader2d}).  Note that $ \{ {\bf x}_p, {\bf y}_p \} $
and $ \{ {\bf x}_e, {\bf y}_e \} $ are planar natural Frenet frames
for the trajectories of the pursuer and evader, respectively.

We model the pursuer and evader as point particles (confined
to the plane), and use natural frames and curvature controls to describe
their motion, because this is a simple model for which we
can derive both physical intuition and concrete control laws.
(Furthermore, although we save the details for a future paper,
this approach generalizes nicely for three-dimensional motion.)
Flying insects and animals (also unmanned aerial vehicles) have
limited maneuverability and must maintain sufficient airspeed to
stay aloft, so treating their motion as constant-speed with
steering control is physically reasonable, at least for
some range of flight conditions.  (Note that the
steering control directly drives the angular velocity of
the particle, and hence is actually an acceleration input.
However, this acceleration is constrained to be perpendicular to the
instantaneous direction of motion, and therefore the speed remains 
unchanged.)

We refer to (\ref{pursuer2d}) and (\ref{evader2d}) as the 
``pursuit-evader system.''  In what follows, we assume that
the pursuer follows a {\it feedback} strategy to drive the system
toward a state of motion
camouflage, and close in on the evader.  The evader, on the other
hand, follows an open-loop strategy.
The analysis we present for the pursuer feedback strategy also
suggests (with a sign change in the control law) how the
evader could use feedback and a motion-camouflage strategy 
to conceal its flight from the pursuer.  Ultimately, it would
be interesting to address the game-theoretic problem in which
both the pursuer and evader follow feedback strategies, so that
the system would truly be a pursuit-evader system.  (What
we address in this work would be more properly described as a
pursuer-pursuee system.  However, we keep the pursuer-evader
terminology, because it sets the stage for 
analyzing the true pursuer-evader system, which we plan to address
in a future paper.)

\subsection{Motion camouflage with respect to the point at infinity}

Motion camouflage with respect to the point at infinity is given by
\begin{equation}
{\bf r}_p = {\bf r}_e + \lambda {\bf r}_\infty,
\end{equation}
where $ {\bf r}_\infty $ is a fixed unit vector and $ \lambda $ is
a time-dependent scalar (see also Section 5 of \cite{glendinning}).

Let 
\begin{equation}
\label{rdefnplanar}
{\bf r} = {\bf r}_p - {\bf r}_e 
\end{equation}
be the vector from the
evader to the pursuer.  We refer to $ {\bf r} $ as the 
``baseline vector,'' and $ |{\bf r}| $ as the ``baseline length.''
We restrict attention to non-collision states, i.e., $ {\bf r} \ne 0 $.
In that case, the component of the pursuer velocity $ \dot{\bf r}_p $
transverse to the base line is
\begin{equation*}
\dot{\bf r}_p - \left( \frac{\bf r}{|{\bf r}|} \cdot \dot{\bf r}_p
 \right) \frac{\bf r}{|{\bf r}|},
\end{equation*}
and similarly, that of the evader is
\begin{equation*}
\dot{\bf r}_e - \left( \frac{\bf r}{|{\bf r}|} \cdot \dot{\bf r}_e
 \right) \frac{\bf r}{|{\bf r}|}.
\end{equation*}
The {\it relative} transverse component is
\begin{eqnarray}
\label{wdefn}
{\bf w} \hspace{-.2cm} & = & \hspace{-.2cm}
\left(\dot{\bf r}_p - \dot{\bf r}_e\right)
 - \left( \frac{\bf r}{|{\bf r}|} \cdot \left(\dot{\bf r}_p - \dot{\bf r}_e
 \right) \right) \frac{\bf r}{|{\bf r}|} \nonumber \\
\hspace{-.2cm} & = & \hspace{-.2cm}
\dot{\bf r} - \left( \frac{\bf r}{|{\bf r}|} \cdot \dot{\bf r}
 \right) \frac{\bf r}{|{\bf r}|}.
\end{eqnarray}

\vspace{.5cm}

\noindent
{\bf Lemma} (Infinitesimal characterization of motion camouflage):
The pursuit-evasion system (\ref{pursuer2d}), (\ref{evader2d})
is in a state of motion camouflage without collision on
an interval iff $ {\bf w} = 0 $ on that interval.

\vspace{.5cm}

\noindent
{\bf Proof}: $ (\Longrightarrow) $ Suppose motion camouflage holds.
Thus
\begin{equation}
{\bf r}(t) = \lambda(t) {\bf r}_\infty, \;\; t \in [0,T].
\end{equation}
Differentiating, $ \dot{\bf r} = \dot{\lambda} {\bf r}_\infty $.  Hence,
\begin{eqnarray}
{\bf w} \hspace{-.2cm} & = & \hspace{-.2cm}
 \dot{\bf r} - \left( \frac{\bf r}{|{\bf r}|} \cdot \dot{\bf r}
 \right) \frac{\bf r}{|{\bf r}|} \nonumber \\
\hspace{-.2cm} & = & \hspace{-.2cm}
\dot{\lambda} {\bf r}_\infty - \left( \frac{\lambda}{|\lambda|} {\bf r}_\infty
 \cdot \dot{\lambda} {\bf r}_\infty \right) \frac{\lambda}{|\lambda|}
 {\bf r}_\infty \nonumber \\
 \hspace{-.2cm} & = & \hspace{-.2cm}
 0  \mbox{ on } [0,T].
\end{eqnarray}
$ (\Longleftarrow) $ Suppose $ {\bf w}=0 $ on $ [0,T] $.  Thus
\begin{equation}
\dot{\bf r} =  \left( \frac{\bf r}{|{\bf r}|} \cdot \dot{\bf r}
 \right) \frac{\bf r}{|{\bf r}|} \triangleq \xi {\bf r},
\end{equation}
so that
\begin{eqnarray}
{\bf r}(t) \hspace{-.2cm} & = & \hspace{-.2cm}
 \exp\left( \int_0^t \xi(\sigma) d\sigma \right) {\bf r}(0)
 \nonumber \\
\hspace{-.2cm} & = & \hspace{-.2cm}
 |{\bf r}(0)| \exp\left( \int_0^t \xi(\sigma) d\sigma \right)
 \frac{{\bf r}(0)}{|{\bf r}(0)|} \nonumber \\
\hspace{-.2cm} & = & \hspace{-.2cm}
 \lambda(t) {\bf r}_\infty,
\end{eqnarray}
where $ {\bf r}_\infty  = {\bf r}(0)/|{\bf r}(0)| $ and
$ \lambda(t) = |{\bf r}(0)| \exp\left( \int_0^t \xi(\sigma) d\sigma \right) $.
$ \Box $

\vspace{.5cm}

It follows from the {\bf Lemma} that the set of all motion
camouflage states constitutes a 5-dimensional smooth manifold
with two connected components, each diffeomorphic to 
$ S^1 \times \mathbb{R} \times SE(2) $ in the 6-dimensional
state space $ SE(2) \times SE(2) $ of the problem.
In practice we are interested in how far the pursuit-evasion system 
is from a state of motion camouflage.  In what follows, we offer a
measure of this.

\subsection{Cost function}

Consider the ratio 
\begin{equation}
\Gamma(t) = \frac{ \frac{d}{dt}{|{\bf r}|}}{\left|\frac{d{\bf r}}{dt}\right|},
\end{equation}
which compares the rate of change of the 
baseline length to the absolute rate of change of the baseline vector.
If the baseline experiences pure lengthening, then the ratio assumes
its maximum value, $ \Gamma(t) = 1 $.
If the baseline experiences pure shortening, then the ratio assumes
its minimum value, $ \Gamma(t) = -1 $.  If the baseline experiences pure
rotation, but remains the same length, then $ \Gamma(t) = 0 $.
Noting that
\begin{equation}
\frac{d}{dt}{|{\bf r}|} = \frac{\bf r}{|{\bf r}|}\cdot \dot{\bf r},
\end{equation}
we see that $ \Gamma(t) $ may alternatively be written as
\begin{equation}
\label{gammadotprod}
\Gamma(t) = \frac{\bf r}{|{\bf r}|} \cdot \frac{\dot{\bf r}}{|\dot{\bf r}|}.
\end{equation}
Thus, $ \Gamma(t) $ is the dot product of two unit vectors: one in the 
direction of $ {\bf r} $, and the other in the direction of $ \dot{\bf r} $.
Note that $ \Gamma $ is well-defined except at $ {\bf r} = 0 $, since
\begin{equation}
\label{dotrbound}
1-\nu \le |\dot{\bf r}| = |{\bf x}_p - \nu {\bf x}_e| \le 1+\nu.
\end{equation} 
For convenience, we define the notation $ {\bf q}^{\perp} $ to represent
the vector $ {\bf q} $ rotated counter-clockwise in the plane by
an angle $ \pi/2 $.  Thus, for example, $ {\bf x}_p^{\perp} = {\bf y}_p $.
The transverse component $ {\bf w} $ of relative velocity, expression
(\ref{wdefn}), then becomes
\begin{eqnarray}
\label{transrelvel2}
{\bf w} \hspace{-.2cm} & = & \hspace{-.2cm}
 \dot{\bf r} - \left(\frac{\bf r}{|{\bf r}|}\cdot \dot{\bf r} \right)
 \frac{\bf r}{|{\bf r}|} = 
\left[ \left(\frac{\bf r}{|{\bf r}|}\right)^{\perp} \cdot \dot{\bf r} \right]
\left(\frac{\bf r}{|{\bf r}|}\right)^{\perp} \nonumber \\
 \hspace{-.2cm} & = & \hspace{-.2cm}
  -  \left(\frac{\bf r}{|{\bf r}|} \cdot \dot{\bf r}^{\perp}\right)
\left(\frac{\bf r}{|{\bf r}|}\right)^{\perp}.
\end{eqnarray}
For convenience, we define $ w $ to be the (signed) magnitude of $ {\bf w} $,
i.e.,
\begin{equation}
\label{transrelvel3}
w = {\bf w} \cdot \left(\frac{\bf r}{|{\bf r}|}\right)^{\perp}
 = -\left(\frac{\bf r}{|{\bf r}|} \cdot \dot{\bf r}^{\perp}\right),
\end{equation}
and refer also to $ w $ as the transverse component of the relative
velocity.
From the orthogonal decomposition
\begin{equation}
\frac{\bf r}{|{\bf r}|}
 =  \left(\frac{\bf r}{|{\bf r}|}\cdot \frac{\dot{\bf r}}{|\dot{\bf r}|}\right)
  \left(\frac{\dot{\bf r}}{|\dot{\bf r}|}\right)
 + \left[\frac{\bf r}{|{\bf r}|}\cdot 
 \left(\frac{\dot{\bf r}}{|\dot{\bf r}|}\right)^\perp \right]
  \left(\frac{\dot{\bf r}}{|\dot{\bf r}|}\right)^\perp,
\end{equation}
it follows that
\begin{equation}
1 =
 \left(\frac{\bf r}{|{\bf r}|}\cdot \frac{\dot{\bf r}}{|\dot{\bf r}|}\right)^2
 + \left[\frac{\bf r}{|{\bf r}|}\cdot 
 \left(\frac{\dot{\bf r}}{|\dot{\bf r}|}\right)^\perp \right]^2
 = \Gamma^2 + \frac{|w|^2}{|\dot{\bf r}|^2}.
\end{equation}
Thus $ (1-\Gamma^2) $ is a measure of the distance from motion camouflage.

\subsection{Feedback law derivation}

Differentiating $ \Gamma $ along trajectories of (\ref{pursuer2d})
and (\ref{evader2d}) gives
\begin{eqnarray}
\label{dotgamma}
\dot{\Gamma} \hspace{-.2cm} & = & \hspace{-.2cm}
\left(\frac{\dot{\bf r} \cdot \dot{\bf r} + {\bf r} \cdot \ddot{\bf r}}
{|{\bf r}||\dot{\bf r}|} \right)
 - \left(\frac{{\bf r}\cdot \dot{\bf r}}{|\dot{\bf r}|} \right)
 \left( \frac{{\bf r} \cdot \dot{\bf r}}{|{\bf r}|^3} \right)
 - \left( \frac{{\bf r}\cdot \dot{\bf r}}{|{\bf r}|} \right)
 \left( \frac{\dot{\bf r} \cdot \ddot{\bf r}}{|\dot{\bf r}|^3} \right)
\nonumber \\
\hspace{-.2cm} & = & \hspace{-.2cm}
\frac{|\dot{\bf r}|}{|{\bf r}|}\left[1 - \left(\frac{\bf r}{|{\bf r}|}
 \cdot \frac{\dot{\bf r}}{|\dot{\bf r}|}\right)^2\right]
 \nonumber \\ & & +\frac{1}{|\dot{\bf r}|}\left[ \frac{\bf r}{|{\bf r}|} -
 \left(\frac{\bf r}{|{\bf r}|}\cdot \frac{\dot{\bf r}}{|\dot{\bf r}|}\right)
 \frac{\dot{\bf r}}{|\dot{\bf r}|} \right] \cdot \ddot{\bf r}.
\end{eqnarray}

From (\ref{rdefnplanar}) we obtain
\begin{equation}
\dot{\bf r}^{\perp} =  {\bf y}_p - \nu {\bf y}_e,
\end{equation}
and
\begin{equation}
\label{ddotrplanar}
\ddot{\bf r} = {\bf y}_p u_p - \nu^2 {\bf y}_e u_e.
\end{equation}
Also,
\begin{eqnarray}
\left[ \frac{\bf r}{|{\bf r}|}
 -  \left(\frac{\bf r}{|{\bf r}|}\cdot \frac{\dot{\bf r}}{|\dot{\bf r}|}\right)
 \frac{\dot{\bf r}}{|\dot{\bf r}|} \right]
 \hspace{-.2cm} & = & \hspace{-.2cm}
 \left[\frac{\bf r}{|{\bf r}|} \cdot
 \left(\frac{\dot{\bf r}}{|\dot{\bf r}|}\right)^{\perp} \right]
 \left(\frac{\dot{\bf r}}{|\dot{\bf r}|}\right)^{\perp} \nonumber \\
 \hspace{-.2cm} & = & \hspace{-.2cm}
\frac{1}{|\dot{\bf r}|^2}
 \left(\frac{\bf r}{|{\bf r}|} \cdot \dot{\bf r}^{\perp} \right)
 \dot{\bf r}^{\perp}.
\end{eqnarray}
Then from (\ref{dotgamma}) we obtain
\begin{eqnarray}
\label{finaldotgammaplanar}
\dot{\Gamma} \hspace{-.2cm} & = & \hspace{-.2cm}
\frac{|\dot{\bf r}|}{|{\bf r}|}\left[1 - \left(\frac{\bf r}{|{\bf r}|}
 \cdot \frac{\dot{\bf r}}{|\dot{\bf r}|}\right)^2\right]
\nonumber \\ & & 
+\frac{1}{|\dot{\bf r}|}\left[ \frac{\bf r}{|{\bf r}|} -
 \left(\frac{\bf r}{|{\bf r}|}\cdot \frac{\dot{\bf r}}{|\dot{\bf r}|}\right)
 \frac{\dot{\bf r}}{|\dot{\bf r}|} \right] \cdot
 \big( {\bf y}_p u_p - \nu^2 {\bf y}_e u_e \big)
\nonumber \\
 \hspace{-.2cm} & = & \hspace{-.2cm}
\frac{|\dot{\bf r}|}{|{\bf r}|}\left[
 \frac{1}{|\dot{\bf r}|^2}
 \left(\frac{\bf r}{|{\bf r}|} \cdot \dot{\bf r}^{\perp} \right)^2 \right]
 \nonumber \\ & &
+\frac{1}{|\dot{\bf r}|}\left[\frac{1}{|\dot{\bf r}|^2}
 \left(\frac{\bf r}{|{\bf r}|} \cdot \dot{\bf r}^{\perp} \right)
 \dot{\bf r}^{\perp}\right]
\cdot \big( {\bf y}_p u_p - \nu^2 {\bf y}_e u_e \big).
\nonumber \\
\end{eqnarray}
Noting that
\begin{equation}
\dot{\bf r}^{\perp} \cdot {\bf y}_p =
\dot{\bf r} \cdot {\bf x}_p = 1 - \nu ({\bf x}_p \cdot {\bf x}_e)
 \ge 1-\nu > 0,
\end{equation}
and
\begin{equation}
\dot{\bf r}^{\perp} \cdot {\bf y}_e =
\dot{\bf r} \cdot {\bf x}_e = ({\bf x}_p \cdot {\bf x}_e) - \nu, 
\end{equation}
we obtain
\begin{eqnarray}
\label{dotgammaupue}
\dot{\Gamma} \hspace{-.2cm} & = & \hspace{-.2cm}
\frac{|\dot{\bf r}|}{|{\bf r}|}\left[
 \frac{1}{|\dot{\bf r}|^2}
 \left(\frac{\bf r}{|{\bf r}|} \cdot \dot{\bf r}^{\perp} \right)^2 \right]
 \nonumber \\ & &
+\frac{1}{|\dot{\bf r}|}\left[\frac{1}{|\dot{\bf r}|^2}
 \left(\frac{\bf r}{|{\bf r}|} \cdot \dot{\bf r}^{\perp} \right) \right]
 \big(1 - \nu ({\bf x}_p \cdot {\bf x}_e) \big) u_p \nonumber \\ & &
 + \frac{1}{|\dot{\bf r}|}\left[\frac{1}{|\dot{\bf r}|^2}
 \left(\frac{\bf r}{|{\bf r}|} \cdot \dot{\bf r}^{\perp} \right) \right]
 \big(\nu - ({\bf x}_p \cdot {\bf x}_e) \big) \nu^2 u_e.
\end{eqnarray}
Suppose that we take
\begin{equation}
\label{planarup}
 u_p = -\mu \left(\frac{\bf r}{|{\bf r}|} \cdot \dot{\bf r}^{\perp} \right)
 + \left[\frac{({\bf x}_p \cdot {\bf x}_e) - \nu }
 { 1-\nu ({\bf x}_p \cdot {\bf x}_e) } \right] \nu^2 u_e,
\end{equation}
where $ \mu > 0 $,
so that the steering control for the pursuer consists of two terms: one
involving the motion of the evader, and one involving
the transverse component of the relative velocity.
Then
\begin{equation}
\label{dotgammaupctrl}
\dot{\Gamma} =  
 - \left[\frac{\mu}{|\dot{\bf r}|}\big(1-\nu ({\bf x}_p \cdot {\bf x}_e)\big)
 - \frac{|\dot{\bf r}|}{|{\bf r}|} \right]
 \left[ \frac{1}{|\dot{\bf r}|}
 \left(\frac{\bf r}{|{\bf r}|} \cdot \dot{\bf r}^{\perp} \right) \right]^2,
\end{equation}
and for any choice of $ \mu > 0 $, there exists $ r_o > 0 $ such that
\begin{equation}
\frac{\mu}{|\dot{\bf r}|}\big(1-\nu ({\bf x}_p \cdot {\bf x}_e)\big)
 - \frac{|\dot{\bf r}|}{|{\bf r}|} > 0,
\end{equation}
for all $ {\bf r} $ such that $ |{\bf r}|> r_o $.  Thus, for control law
(\ref{planarup}), 
\begin{equation}
\dot{\Gamma} \le 0, \;\; \forall |{\bf r}|> r_o.
\end{equation}

\section{The high-gain limit}

Control law (\ref{planarup}) has the nice property that for any value
of the gain $ \mu >0 $, there is a disc of radius $ r_o $ (depending on
$ \mu $) such that $ \dot{\Gamma} \le 0 $ outside the disc.
However, the problem with (\ref{planarup}) is that the pursuer needs
to know (i.e., sense and estimate) the evader's steering program $ u_e $.
Here we show that by taking $ \mu $ sufficiently large, motion camouflage
can be achieved (in a sense we will make precise) using a control
law depending only on the transverse relative velocity:
\begin{equation}
\label{planarupgain}
u_p = -\mu \left(\frac{\bf r}{|{\bf r}|} \cdot \dot{\bf r}^{\perp} \right),
\end{equation}
in place of (\ref{planarup}), provided $ |u_e| $ is bounded.
Comparing (\ref{planarupgain}) to (\ref{transrelvel3}), we
see that, indeed, $ u_p $ is proportional to the signed length of
the relative transverse velocity vector.
We will designate this as the {\it motion camouflage proportional guidance}
(MCPG) law for future reference (see Section V below).

As is further discussed in Section V, (\ref{planarupgain})
requires range information as well as pure optical flow sensing.
However, the range information can be coarse, since range 
errors (within appropriate bounds) have the same effect in
(\ref{planarupgain}) as gain variations.  
We say that (\ref{planarupgain}) is {\it biologically
plausible} because the only critical sensor measurement required is
optical flow sensing.
Optical flow sensing does not yield the relative transverse
velocity directly, but rather the angular speed of the 
image of the evader across the pursuer's eye.  In fact,
it is the sign of the optical flow that is most critical
to measure correctly, since errors in the magnitude of the
optical flow, like range errors, only serve to modulate the gain
in (\ref{planarupgain}).

For biological systems, the capabilities of the sensors {\it vis-a-vis}
the sensing requirements for implementing (\ref{planarupgain})
constrain the range of conditions for which (\ref{planarupgain})
represents a feasible control strategy.  In the high-gain limit
we focus on below, sensor noise (which is amplified by the high gain)
would be expected to have significant impact.  However, to illustrate
the essential behavior, here we neglect both sensor limitations and noise.

\subsection{Bounds and estimates}

Let us consider control law (\ref{planarupgain}), and the resulting behavior
of $ \Gamma $ as a function of time.
From (\ref{dotgammaupue}), we obtain the inequality
\begin{eqnarray}
\label{dotgammabound}
\dot{\Gamma} \hspace{-.2cm} & = & \hspace{-.2cm}
 - \left[\frac{\mu}{|\dot{\bf r}|}\big(1-\nu ({\bf x}_p \cdot {\bf x}_e)\big)
 - \frac{|\dot{\bf r}|}{|{\bf r}|} \right]
 \left[ \frac{1}{|\dot{\bf r}|}
 \left(\frac{\bf r}{|{\bf r}|} \cdot \dot{\bf r}^{\perp} \right) \right]^2
 \nonumber \\ & & 
 + \frac{1}{|\dot{\bf r}|}\left[\frac{1}{|\dot{\bf r}|^2}
 \left(\frac{\bf r}{|{\bf r}|} \cdot \dot{\bf r}^{\perp} \right) \right]
 \big(\nu - ({\bf x}_p \cdot {\bf x}_e) \big) \nu^2 u_e
 \nonumber \\
\hspace{-.2cm} & \le & \hspace{-.2cm}
 - \left(1-\Gamma^2 \right)
 \left[\frac{\mu}{|\dot{\bf r}|}\big(1-\nu ({\bf x}_p \cdot {\bf x}_e)\big)
 - \frac{|\dot{\bf r}|}{|{\bf r}|} \right] \nonumber \\ & &
 + \frac{1}{|\dot{\bf r}|^2}\sqrt{1-\Gamma^2}
 \Big| \big(\nu - ({\bf x}_p \cdot {\bf x}_e) \big) \nu^2 u_e \Big|
 \nonumber \\
\hspace{-.2cm} & \le & \hspace{-.2cm}
 - \left(1-\Gamma^2 \right)
 \left[\frac{\mu}{|\dot{\bf r}|}(1-\nu)
 - \frac{|\dot{\bf r}|}{|{\bf r}|} \right] \nonumber \\ & &
 + \left(\sqrt{1-\Gamma^2}\right)\frac{\nu^2(1+\nu)(\max |u_e|)}
 {|\dot{\bf r}|^2}
 \nonumber \\
\hspace{-.2cm} & \le & \hspace{-.2cm}
  - \left(1-\Gamma^2 \right)
 \left[\mu\left(\frac{1-\nu}{1+\nu}\right)
 - \frac{1+\nu}{|{\bf r}|} \right] \nonumber \\ & &
 + \left(\sqrt{1-\Gamma^2}\right)\left[\frac{\nu^2(1+\nu)(\max |u_e|)}
 {(1-\nu)^2}\right],
\end{eqnarray}
where we have used (\ref{dotrbound}).
For convenience, we define the constant $ c_1 > 0 $ as
\begin{equation}
c_1 = \frac{\nu^2(1+\nu)(\max |u_e|)}{(1-\nu)^2}.
\end{equation}

For any $ \mu > 0 $, we can define
$ r_o > 0 $ and $ c_o > 0 $ such that
\begin{equation}
\label{mudecomp}
\mu  = \left(\frac{1+\nu}{1-\nu}\right) \left(
 \frac{1+\nu}{r_o} + c_o \right)
\end{equation}
(and it is clear that many such choices of $ r_o $ and $ c_o $ exist).
Note that (\ref{mudecomp}) implies
\begin{equation}
\mu \ge \left(\frac{1+\nu}{1-\nu}\right) \left(
 \frac{1+\nu}{|{\bf r}|} + c_o \right), \;\; \forall |{\bf r}| \ge r_o.
\end{equation}
Thus, for $ |{\bf r}| \ge r_o $, (\ref{dotgammabound}) becomes
\begin{eqnarray}
\dot{\Gamma} \hspace{-.2cm} & \le & \hspace{-.25cm}
   - \left(1-\Gamma^2 \right) \hspace{-.1cm}
 \left[ \hspace{-.05cm} \left(\frac{1+\nu}{1-\nu}\right) \hspace{-.1cm} \left(
 \frac{1+\nu}{|{\bf r}|} + c_o \right) \hspace{-.1cm}
 \left(\frac{1-\nu}{1+\nu}\right)
 - \frac{1+\nu}{|{\bf r}|} \hspace{-.05cm} \right] \nonumber \\ & &
 + \left(\sqrt{1-\Gamma^2}\right)c_1
  \nonumber \\
 \hspace{-.2cm} & = & \hspace{-.2cm}
   - \left(1-\Gamma^2 \right) c_o + \left(\sqrt{1-\Gamma^2} \right)c_1.
\end{eqnarray}
Suppose that given $ 0<\epsilon << 1$,
we take $ c_o \ge 2c_1/\sqrt{\epsilon} $.  Then for $ (1-\Gamma^2) > \epsilon $,
\begin{eqnarray}
\label{dotgammabound2}
\dot{\Gamma}  \hspace{-.2cm} & \le & \hspace{-.25cm} 
 - \left(1-\Gamma^2 \right) c_o + \left(\sqrt{1-\Gamma^2} \right)c_1
 \nonumber \\
 \hspace{-.2cm} & = & \hspace{-.2cm}
  - \left(1-\Gamma^2 \right) \left( c_o
 - \frac{c_1}{\sqrt{1-\Gamma^2}} \right) \nonumber \\
  \hspace{-.2cm} & \le & \hspace{-.25cm}
  - \left(1-\Gamma^2 \right) \left( c_o
 - \frac{c_1}{\sqrt{\epsilon}}\right) \nonumber \\
 \hspace{-.2cm} & = & \hspace{-.2cm}
 -  \left(1-\Gamma^2 \right) c_2,
\end{eqnarray} 
where
\begin{equation}
\label{c2defn}
c_2 = c_o - \frac{c_1}{\sqrt{\epsilon}} > 0.
\end{equation}

\vspace{.5cm}

\noindent
{\bf Remark}: There are two possibilities for 
\begin{equation}
\label{oneminusgammasqreps}
 (1-\Gamma^2) \le \epsilon.
\end{equation} 
The state we seek to drive the system toward has $ \Gamma \approx -1 $;
however, (\ref{oneminusgammasqreps}) can also be satisfied for 
$ \Gamma \approx 1 $.  (Recall that $ -1 \le \Gamma \le 1 $.)
There is always a set of initial conditions
such that (\ref{oneminusgammasqreps}) is satisfied with $ \Gamma \approx 1 $.
We can address this issue as follows: let $ \epsilon_o > 0 $ denote how
close to $ -1 $ we wish to drive $ \Gamma $, and let $ \Gamma_0 = \Gamma(0) $
denote the initial value of $ \Gamma $.  Take
\begin{equation}
\epsilon = \min(\epsilon_o,1-\Gamma_0^2),
\end{equation}
so that (\ref{dotgammabound2}) with (\ref{c2defn})
applies from time $ t = 0 $. $ \Box $

\vspace{.5cm}

From (\ref{dotgammabound2}), we can write
\begin{equation}
\frac{d\Gamma}{1-\Gamma^2} \le - c_2 dt,
\end{equation}
which, integrating both sides, leads to
\begin{equation}
\int_{\Gamma_0}^{\Gamma} \frac{d\tilde{\Gamma}}{1-\tilde{\Gamma}^2} \le 
 -c_2 \int_0^t d\tilde{t} = -c_2 t,
\end{equation}
where $ \Gamma_0 = \Gamma(t=0) $.
Noting that
\begin{equation}
\int_{\Gamma_0}^{\Gamma}  \frac{d\tilde{\Gamma}}{1-\tilde{\Gamma}^2}
 = \int_{\Gamma_0}^{\Gamma} d(\tanh^{-1} \tilde{\Gamma}) 
 = \tanh^{-1} \Gamma - \tanh^{-1} \Gamma_0,
\end{equation}
we see that for $ |{\bf r}| \ge r_o $, (\ref{dotgammabound2}) implies
\begin{equation}
\label{gammabound}
\Gamma(t) \le \tanh\left(\tanh^{-1}\Gamma_0 - c_2 t \right),
\end{equation}
where we have used the fact that $ \tanh^{-1}(\cdot) $ is
a monotone increasing function.

Now we consider estimating how long $ |{\bf r}| \ge r_o $, 
which in turn determines how large $ t $ can become in inequality
(\ref{gammabound}), and hence how close to $-1$ will $ \Gamma(t) $ be 
driven.  From (\ref{gammadotprod}) we have
\begin{equation}
\frac{d}{dt}|{\bf r}| = \Gamma(t) |\dot{\bf r}|,
\end{equation}
which from (\ref{dotrbound}) and $ |\Gamma(t)| \le 1 $, $ \forall t $,
implies 
\begin{equation}
\label{dotnormrbound}
\frac{d}{dt}|{\bf r}| \ge -|\Gamma(t)|(1+\nu) \ge - (1+\nu).
\end{equation}
From (\ref{dotnormrbound}), we conclude that
\begin{equation}
|{\bf r}(t)| \ge |{\bf r}(0)| - (1+\nu)t, \;\; \forall t \ge 0,
\end{equation}
and, more to the point,
\begin{equation}
\label{dotnormrbound2}
|{\bf r}(t)| \ge r_o, \;\; \forall t \le \frac{|{\bf r}(0)|-r_o}{1+\nu}.
\end{equation}
For (\ref{dotnormrbound2}) to be meaningful for the problem at hand,
we assume that $ |{\bf r}(0)| > r_o $.  Then defining
\begin{equation}
\label{bigtdefn}
T = \frac{|{\bf r}(0)|-r_o}{1+\nu} > 0 
\end{equation}
to be the minimum interval of time over which we can guarantee
that $ \dot{\Gamma} \le 0 $, we conclude that
\begin{equation}
\label{gammafinal}
\Gamma(T)  \le  
  \tanh \left( \tanh^{-1} \Gamma_0 - c_2 T \right). 
\end{equation}

From (\ref{gammafinal}), we see that by choosing $ c_2 $ sufficiently
large (which can be accomplished by choosing $ c_o \ge 2c_1/\sqrt{\epsilon} $
sufficiently large), we can force $ \Gamma(T) \le -1+\epsilon $.
Noting that
\begin{equation}
\tanh(x) 
 \le -1+\epsilon
 \Longleftrightarrow x \le \frac{1}{2} \ln 
 \left(\frac{\epsilon}{2 -\epsilon} \right),
\end{equation}
for $ 0 < \epsilon << 1 $,
we see that 
\begin{equation}
\Gamma(T) \le -1 + \epsilon
  \Longleftrightarrow
\tanh^{-1} \Gamma_0 - c_2 T \le  \frac{1}{2} \ln 
 \left(\frac{\epsilon}{2 -\epsilon} \right).
\end{equation}
Thus, if $ c_o \ge 2c_1/\sqrt{\epsilon} $ is taken to be sufficiently
large that
\begin{equation}
\label{c2bound}
c_2 \ge (1+\nu)\frac{\tanh^{-1} \Gamma_0 - \frac{1}{2} \ln
 \left(\frac{\epsilon}{2 -\epsilon} \right)}
 { |{\bf r}(0)| - r_o },
\end{equation}
then we are guaranteed (under the conditions mentioned in the 
above calculations) to achieve $ \Gamma(t_1) \le -1 + \epsilon $
at some finite time $ t_1 \le T $.  

\subsection{Statement of result}

\noindent
{\bf Definition}: Given the system
(\ref{pursuer2d}), (\ref{evader2d}) with $ \Gamma $ defined by
(\ref{gammadotprod}), we say that ``motion camouflage is accessible 
in finite time'' if for any $ \epsilon > 0 $ there exists a 
time $ t_1 > 0 $ such that $ \Gamma(t_1) \le -1 + \epsilon $.
$ \Box $

\vspace{.5cm}

\noindent
{\bf Proposition}: Consider the system
(\ref{pursuer2d}), (\ref{evader2d}) with $ \Gamma $ defined by
(\ref{gammadotprod}), and control law (\ref{planarupgain}),
with the following hypotheses:
\begin{itemize}
\item[$ \hspace{-.05cm} (A1) \hspace{-.05cm} $] $ 0 < \nu < 1 $
 (and $ \nu $ is constant),
\item[$ \hspace{-.05cm} (A2) \hspace{-.05cm} $] 
 $ u_e $ is continuous and $ |u_e| $ is bounded,
\item[$ \hspace{-.05cm} (A3) \hspace{-.05cm} $] 
 $ \Gamma_0 = \Gamma(0) < 1 $, and
\item[$ \hspace{-.05cm} (A4) \hspace{-.05cm} $] $ |{\bf r}(0)| > 0 $.
\end{itemize}
Motion camouflage is accessible in finite time
using high-gain feedback (i.e., by choosing $ \mu > 0 $ sufficiently
large).

\vspace{.5cm}

\noindent
{\bf Proof}: Choose $ r_o > 0 $ such that $ r_o < |{\bf r}(0)| $.
Choose $ c_2 > 0 $ sufficiently large so as to satisfy (\ref{c2bound}),
and choose $ c_o $ accordingly to ensure that (\ref{dotgammabound2}) holds
for $ \Gamma > -1 + \epsilon $.  Then defining $ \mu $ according to
(\ref{mudecomp}) ensures that $ \Gamma(T) \le -1 + \epsilon $,
where $ T > 0 $ is defined by (\ref{bigtdefn}). $ \Box $

\vspace{.5cm}

\noindent
{\bf Remark}: Assumption $ (A1) $ above can be generalized to 
$ 0 \le \nu < 1 $.
(The $ \nu = 0 $ case corresponds to a stationary ``evader,''  so that
the natural Frenet frame (\ref{evader2d}) and steering control $ u_e $
for the evader are not defined.)  $ \Box $

\section{Simulation results}

\begin{figure}
\hspace{.5cm}
\epsfxsize=7cm
\epsfbox{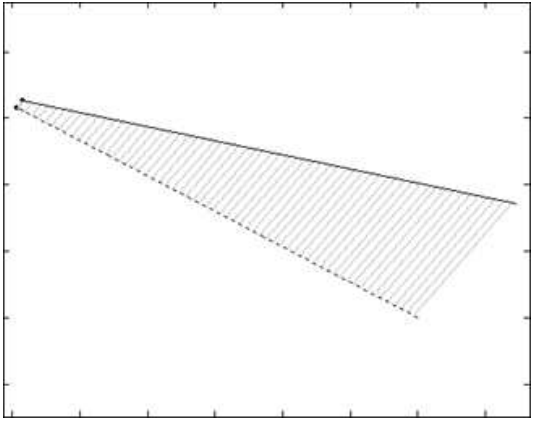}
\caption{\label{straight_traj} Straight-line evader trajectory 
(dashed dark line), and
the corresponding pursuer trajectory (solid dark line) evolving
according to (\ref{pursuer2d}) with control given by
(\ref{planarupgain}).}
\end{figure}

The following simulation results
illustrate the behavior of the pursuit-evasion system (\ref{pursuer2d}),
(\ref{evader2d}), under the control law (\ref{planarupgain}) for the
pursuer and various open-loop controls for the evader.
The simulations also confirm the analytical results presented above.
Figure \ref{straight_traj} shows the behavior of the system for the 
simplest evader behavior, $ u_e = 0 $, which corresponds to 
straight-line motion.  Because control law
(\ref{planarupgain}) is the same as (\ref{planarup}) when $ u_e = 0 $,
$ \Gamma $ tends monotonically toward $ -1 $ (for the initial conditions
and choice of gain $ \mu $ used in the simulation shown).
In figure \ref{straight_traj}, as in the subsequent figures showing
pursuer and evader trajectories, the solid light lines connect the
pursuer and evader positions at evenly-spaced time instants.
For a pursuit-evasion system in a state of motion camouflage, these lines 
would all be parallel to one another.  Also, each simulation is
run for finite time, at the end of which the pursuer and evader
are in close proximity.  (The ratio of speeds is $ \nu = 0.9 $ in all of
the simulations shown.)

\begin{figure}
\hspace{.5cm}
\epsfxsize=7cm
\epsfbox{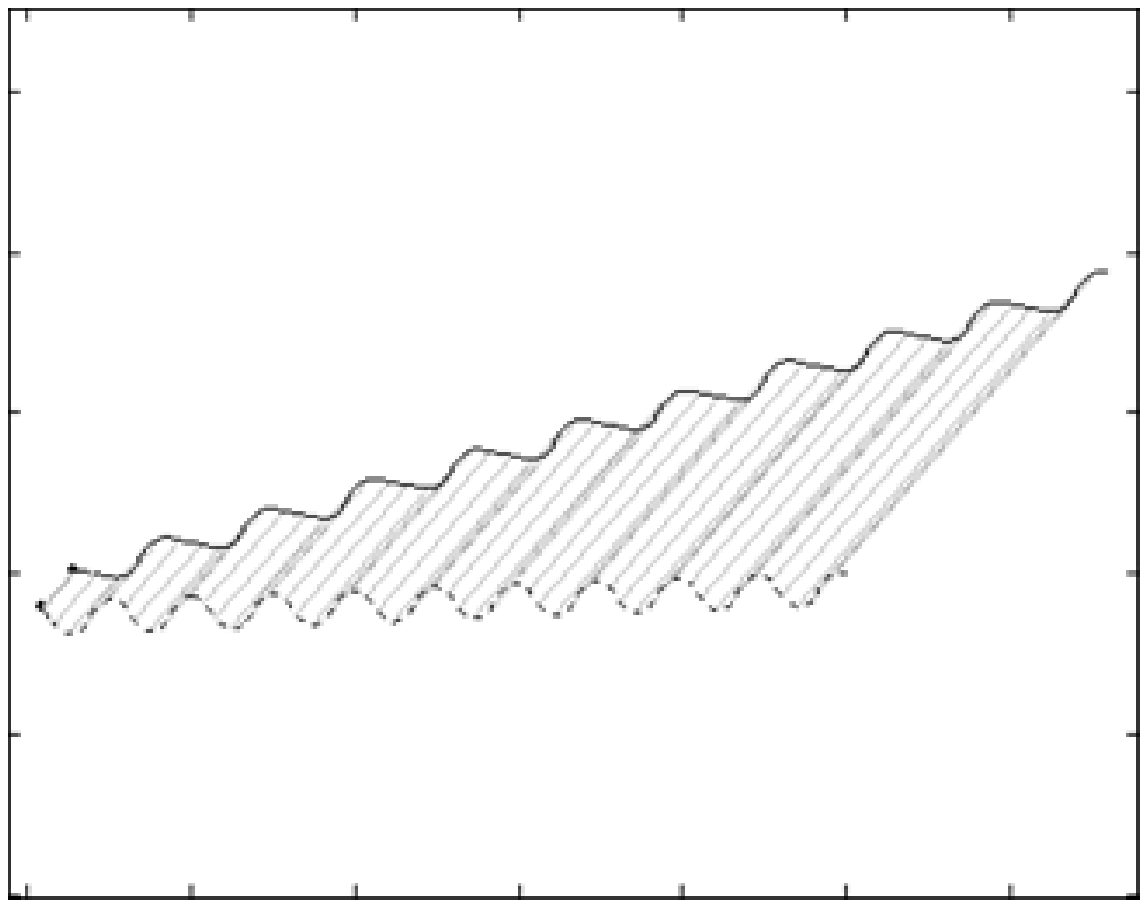}
\caption{\label{sine_traj} Evader trajectory with
sinusoidally varying steering input (dashed dark line), and
the corresponding pursuer trajectory (solid dark line) evolving
according to (\ref{pursuer2d}) with control given by
(\ref{planarupgain}).}
\end{figure}

\begin{figure}
\hspace{.5cm}
\epsfxsize=7cm
\epsfbox{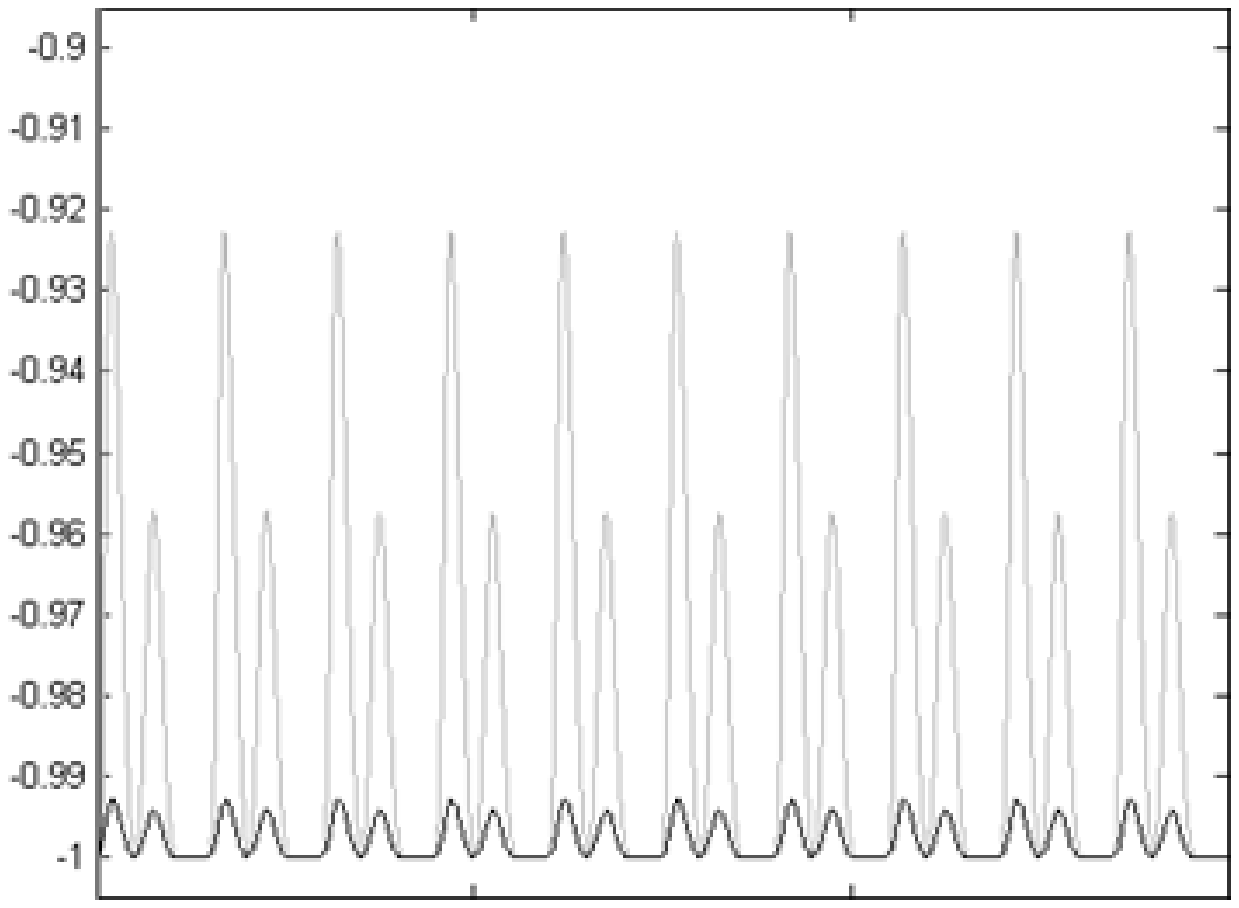}
\caption{\label{sine_gamma} The cost function $ \Gamma(t) $ given
by (\ref{gammadotprod}), plotted as a function of time, for the
pursuit illustrated in figure \ref{sine_traj}.  The two traces
correspond to different values of gain $ \mu $: the value of $ \mu $
is three times as large for the dark trace as for the light trace.
(The trajectories corresponding to the two different gains are
qualitatively similar; figure \ref{sine_traj} actually corresponds to
the lower value of $ \mu $.) }
\end{figure}

\begin{figure}
\hspace{.5cm}
\epsfxsize=7cm
\epsfbox{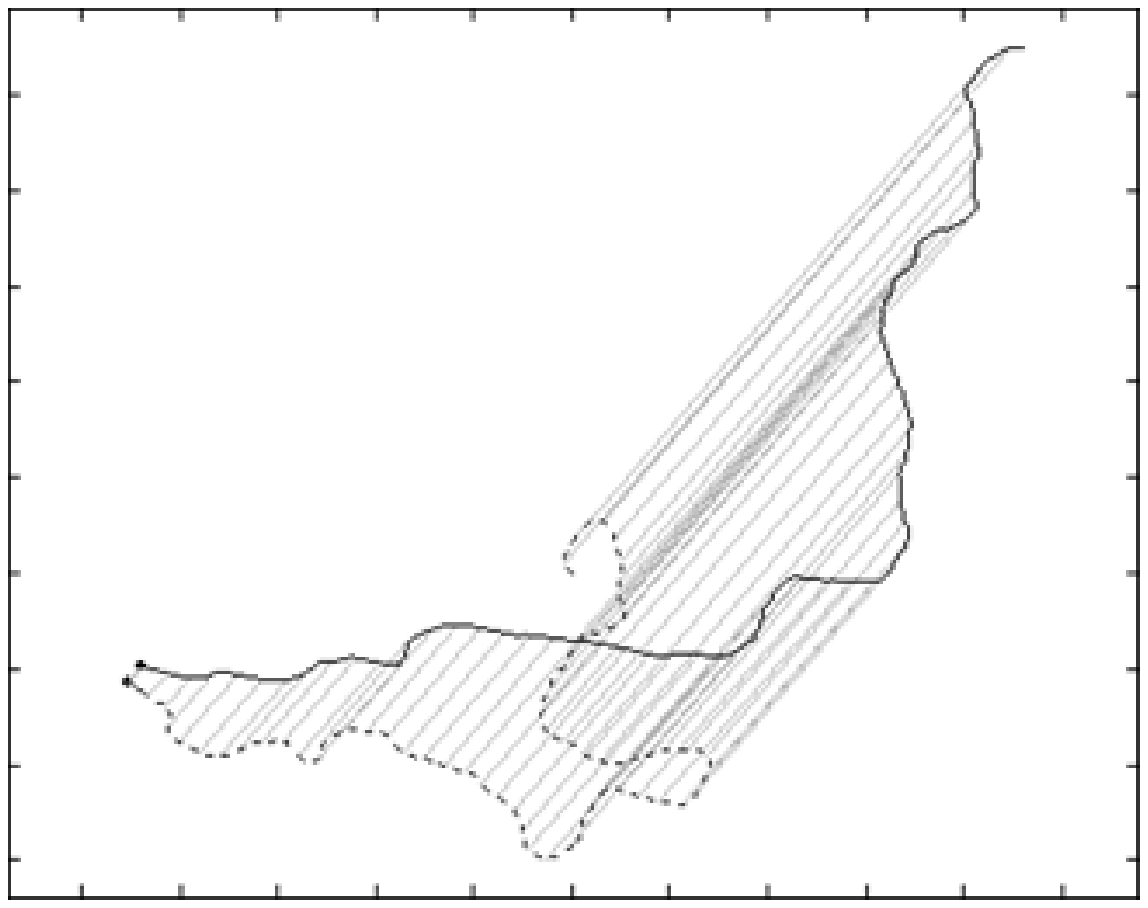}
\caption{\label{random_traj} Evader trajectory with
randomly varying steering input (dashed dark line), and
the corresponding pursuer trajectory (solid dark line) evolving 
according to (\ref{pursuer2d}) with control given by
(\ref{planarupgain}).}
\end{figure}

\begin{figure}
\hspace{.5cm}
\epsfxsize=7cm
\epsfbox{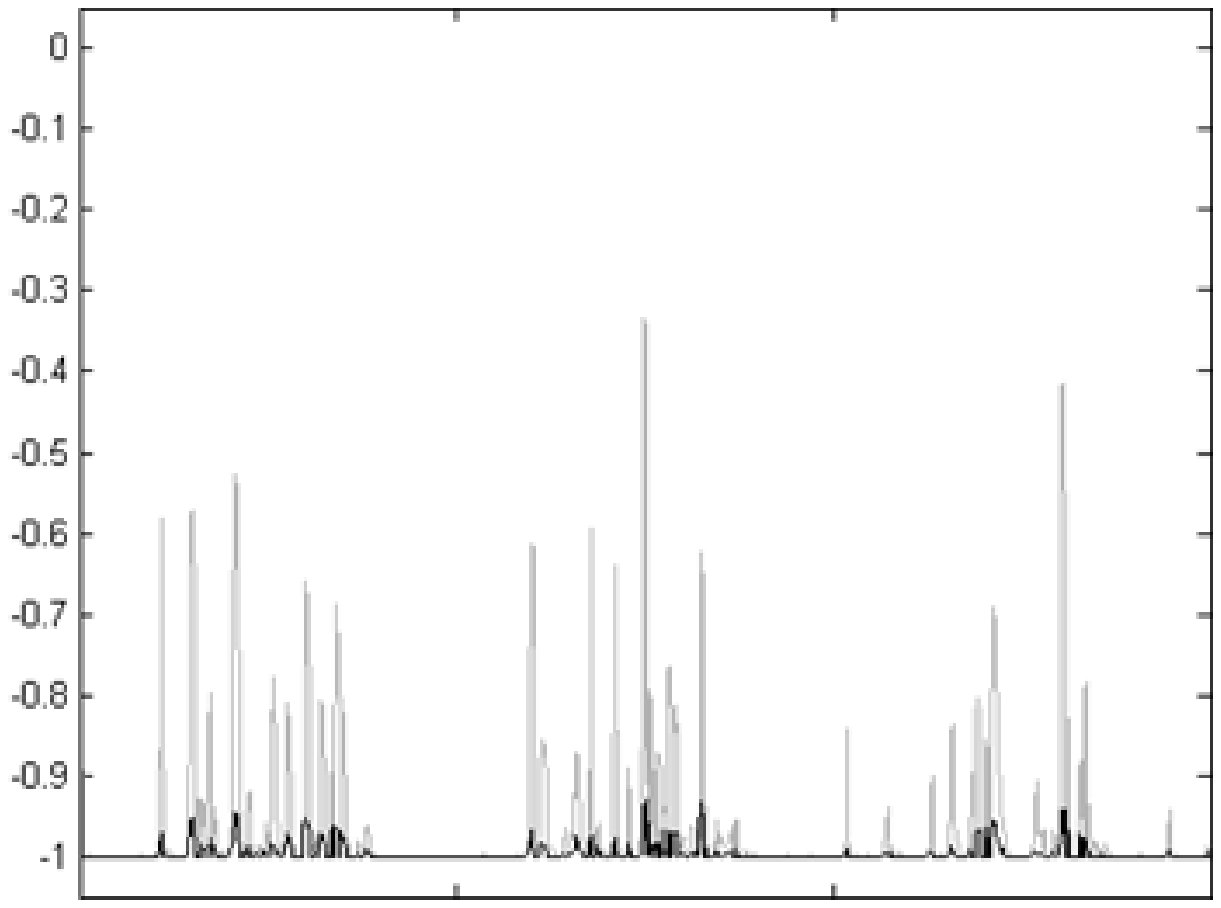}
\caption{\label{random_gamma} The cost function $ \Gamma(t) $ given
by (\ref{gammadotprod}), plotted as a function of time, for the
pursuit illustrated in figure \ref{random_traj}.  The two traces 
correspond to different values of gain $ \mu $: the value of $ \mu $
is three times as large for the dark trace as for the light trace.
(The trajectories corresponding to the two different gains are
qualitatively similar; figure \ref{random_traj} actually corresponds to
the lower value of $ \mu $.) }
\end{figure}

\begin{figure}
\hspace{.5cm}
\epsfxsize=7cm
\epsfbox{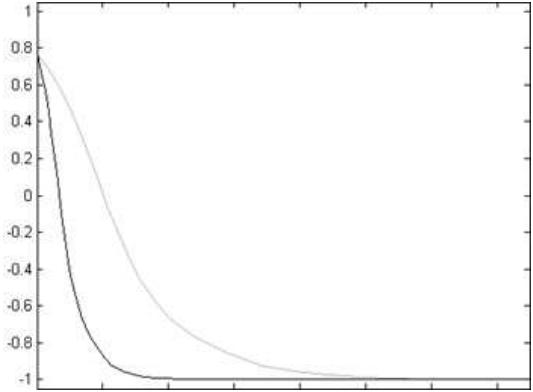}
\caption{\label{random_gamma_initial} The cost function $ \Gamma(t) $ given
by (\ref{gammadotprod}), plotted as a function of time, for 
1/200th of the time interval of figure \ref{random_gamma} for the
pursuit illustrated in figure \ref{random_traj}.  The two traces
correspond to different values of gain $ \mu $: the value of $ \mu $
is three times as large for the dark trace as for the light trace.
(Because similar initial conditions were used, the 
expanded-time-scale plot of $ \Gamma(t) $ corresponding to figures
\ref{sine_traj} and \ref{sine_gamma} is very similar to
figure \ref{random_gamma_initial}.)}
\end{figure}

\begin{figure}
\hspace{.5cm}
\epsfxsize=7cm
\epsfbox{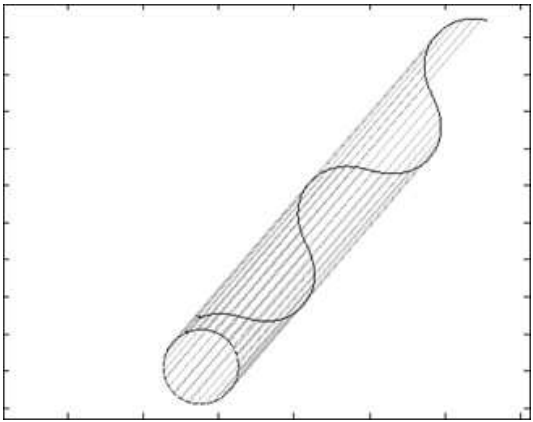}
\caption{\label{circle_traj} Evader trajectory with
constant steering input (circular trajectory), and
the corresponding pursuer trajectory (solid dark line) evolving 
according to (\ref{pursuer2d}) with control given by
(\ref{planarupgain}).}
\end{figure}

Figure \ref{sine_traj} illustrates the behavior of the
pursuer for a sinsusoidally-varying steering control $ u_e $ of the
evader, and figure \ref{sine_gamma} shows the corresponding
behavior of $ \Gamma(t) $.  In figure \ref{sine_gamma}, 
increasing the value of the feedback gain $ \mu $ by a factor of three
is observed to decrease the peak difference between $ \Gamma $ and
$ -1 $ by a factor of about $ 3^2 = 9 $.  This is consistent with
the calculations in the proof of the {\bf Proposition}.
Figure \ref{random_traj} illustrates the behavior of the
pursuer for a randomly-varying steering control $ u_e $ of the
evader, and figures \ref{random_gamma} and \ref{random_gamma_initial}
show the corresponding behavior of $ \Gamma(t) $.  Similarly to figure
\ref{sine_gamma}, figure \ref{random_gamma} shows that increasing
the feedback gain $ \mu $ by a factor of three decreases the peak
difference between $ \Gamma $ and $ -1 $ by a factor of about $ 3^2 = 9 $.
Figure \ref{random_gamma_initial} shows the initial transient in
$ \Gamma(t) $ for $ t $ small.  As would be expected, increasing the
gain $ \mu $ increases the convergence rate.  (The time axes for
figures \ref{random_gamma} and \ref{random_gamma_initial} differ by
a factor of 200, which is why the initial transient cannot be seen
in figure \ref{random_gamma}.)
Finally, figure \ref{circle_traj} illustrates the behavior of the pursuer
for a constant steering control $ u_e $, resulting in circling motion
by the evader.

\section{Connections to missile guidance}

There is a vast literature on the subject of missile guidance in which
the problem of pursuit of an (evasively) maneuvering target by a 
tactical missile is of central interest.  A particular class of
feedback laws, known as {\it pure proportional navigation guidance}
(PPNG) occupies a prominent place \cite{shneydor98}.  For planar
missile-target engagements, the PPNG law determining the steering
control for the missile/pursuer is
\begin{equation}
u^{\mathit PPNG} = N \dot{\lambda},
\end{equation}
where $ \dot{\lambda} $ denotes the rate of rotation (in the plane)
of the line-of-sight (LOS) vector from the pursuer to the evader. 
Here the gain $ N $ is a dimensionless positive constant known as the
navigation constant.  Notice that our motion camouflage guidance law
(MCPG) given by (\ref{planarupgain}) has a gain $ \mu $ which has
the dimensions of $ [\mathit{LENGTH}]^{-1} $.  Also, it is easy to
see that
\begin{equation}
\dot{\lambda} = \frac{w}{|{\bf r}|}
 = -\frac{1}{|{\bf r}|}
 \left( \frac{\bf r}{|{\bf r}|} \cdot \dot{\bf r}^\perp \right).
\end{equation}
So, to make a proper comparison we let $ r_o $ as in Section III be 
a length scale for the problem and define the dimensionless gain
\begin{equation}
N^{\mathit MCPG} = \mu r_o.
\end{equation}
Thus, our MPCG law takes the form
\begin{equation}
u^{\mathit MCPG} = N^{\mathit MCPG} \frac{|{\bf r}|}{r_o} \dot{\lambda}.
\end{equation}
It follows that motion camouflage uses range information to support a
high gain in the initial phase of the engagement, ramping down to a lower
value in the terminal phase $ (|{\bf r}| \approx r_o ) $.  In nature
this extra freedom of gain control is particularly relevant for
echolocating bats (see \cite{ghkm05}), which have remarkable ranging
ability.

Analysis of the performance of the PPNG law is carried out in 
\cite{ha_hur_ko_song90,oh_ha99}, using arguments similar to ours
(although our sufficient conditions appear to be weaker).  While
motion camouflage as a strategy is discussed in \cite{shneydor98},
under ``parallel navigation,'' to the best of our knowledge, the
current work is the first to present and analyze a feedback law
for motion camouflage.

\section{Directions for further work}

In work under preparation, 
we have generalized the analysis to the three-dimensional 
setting, and to planar motion camouflage with respect to a finite point.
The three-dimensional analysis is made possible by the use of natural
Frenet frames, analogously to the three-dimensional 
unit-speed particle interaction laws described in \cite{cdc05}.

Because we are able to treat the motion camouflage problem within
the same framework as our earlier formation control and 
obstacle-avoidance work \cite{scl02}-\cite{cdc05}, we would like
to understand how teams of vehicles can make use of motion camouflage,
and whether we can determine the convergence behavior of such systems.
Various biologically-inspired scenarios for motion camouflage with teams  
have been described in \cite{andersonteam}.   Considering additional 
military applications without biological analogs, there are thus a variety 
of team motion camouflage problems to study.

\section{Acknowledgements}

The authors would like to thank M.V. Srinivasan of the Research School
of Biological Sciences at the Australian National University for
valuable discussions and helpful comments on an earlier draft of this paper.

This research was supported in part by the Naval Research Laboratory under
Grants No.~N00173-02-1G002, N00173-03-1G001, N00173-03-1G019, and
N00173-04-1G014; by the 
Air Force Office of Scientific Research under AFOSR Grants
No.~F49620-01-0415 and FA95500410130; by the Army Research Office
under ODDR\&E MURI01 Program Grant No.~DAAD19-01-1-0465 to the Center for
Communicating Networked Control Systems (through Boston University);
and by NIH-NIBIB grant
1 R01 EB004750-01, as part of the NSF/NIH Collaborative
Research in Computational Neuroscience Program.

\end{document}